\documentclass{elsarticle}

\newtheorem{theorem}{Theorem}[section]

\newtheorem{corollary}[theorem]{Corollary}

\newtheorem{lemma}[theorem]{Lemma}

\newproof{proof}{Proof}

\def\NN{\hbox{\sf I\kern-.13em\hbox{N}}}
\def\RR{\hbox{\sf I\kern-.14em\hbox{R}}}
\def\CC{\hbox{\sf C\kern -.48emC}}
\def\QQ{\hbox{\sf C\kern -.48emQ}}
\def\Cc{\hbox{\sf C\kern -.47em {\raise .48ex \hbox{$\scriptscriptstyle |$}}
   \kern-.5em {\raise .48ex \hbox{$\scriptscriptstyle |$}} }}
\def\Qq{\hbox{\sf Q\kern -.57em {\raise .48ex \hbox{$\scriptscriptstyle |$}}
   \kern-.55em {\raise .48ex \hbox{$\scriptscriptstyle |$}} }}
\def\narr{\advance\leftskip by\parindent}
\def\narrr{\advance\leftskip by-\parindent}

\newcommand{\be}{\begin{equation}}
\newcommand{\ee}{\end{equation}}
\newcommand{\BX}{{\mathcal B}(X)}
\newcommand{\cC}{{\mathcal C}}
\newcommand{\hcC}{\hat{\mathcal C}}

\newcommand{\cM}{{\mathcal M}}
\newcommand{\cN}{{\mathcal N}}

\newcommand{\cS}{{\mathcal S}}
\newcommand{\hT}{\hat{T}}
\newcommand{\cP}{{\mathcal P}}
\newcommand{\cR}{{\mathcal R}}

\newcommand{\rank}{{\rm rank \,}}
\newcommand{\lin}{{\rm lin \,}}

\begin{document}

\title{Invariant subspaces for operator semigroups  
with commutators of rank at most one\tnoteref{t1,t2}}
\tnotetext[t1]{Please cite this article in press as: R. Drnov\v sek, Invariant subspaces for operator semigroups with commutators of
rank at most one, J. Funct. Anal. (2009), doi:10.1016/j.jfa.2009.03.010}
\tnotetext[t2]{This work was supported
by the Slovenian Research Agency.}

\author[rd]{Roman Drnov\v{s}ek}
\ead{roman.drnovsek@fmf.uni-lj.si}
\address[rd]{Department of Mathematics, Faculty of Mathematics and Physics \\
Institute of Mathematics, Physics and Mechanics \\ 
University of Ljubljana \\
Jadranska 19, SI-1000 Ljubljana, Slovenia}


\begin{abstract}
Let $X$ be a complex Banach space of dimension at least $2$, 
and let $\cS$ be a multiplicative semigroup of operators on $X$ such that 
the rank of $ST - TS$ is at most $1$ for all $\{S, T\} \subset \cS$.
We prove that $\cS$ has a non-trivial invariant subspace provided it is 
not commutative. As a consequence we show that $\cS$ is triangularizable
if it consists of polynomially compact operators.
This generalizes results from \cite{RR97} and \cite{CD98}.
\end{abstract}
 
\maketitle

\section{Introduction}

Throughout the paper, let $X$ be a complex Banach space of dimension
at least $2$. A {\it subspace} of $X$ means a closed linear manifold of $X$.
{\it Trivial} subspaces of $X$ are $\{0\}$ and $X$.
The dual space of $X$ is denoted by $X^*$. 
By an {\it operator} on $X$ we mean a bounded linear 
transformation from $X$ into itself. By $I$ we denote the identity operator.
The Banach algebra of all operators on $X$ is denoted by $\BX$.
We denote by $T^*$ the adjoint operator of $T \in \BX$. 
The notation $[S,T]$ is used as an abbreviation for the commutator $S T - T S$,
where $S, T \in \BX$.  
Given $y \in X$ and $\phi \in X^*$, the rank-one operator $y \otimes \phi$ on $X$ is defined 
by $(y \otimes \phi) x = \phi(x) y$.  

A subspace $\cM$ of $X$ is {\it invariant} under an operator $T \in \BX$ 
whenever $T (\cM) \subseteq \cM$. 
Let ${\mathcal C}$ be a collection of operators in $\BX$.
A subspace $\cM$ of $X$ is {\it invariant} under ${\mathcal C}$  
if $\cM$ is invariant under every $T \in \cC$. 
If, in addition, the subspace $\cM$ is invariant under every $S \in \BX$ 
that commutes with all operators of $\cC$, $\cM$ is said to be {\it hyperinvariant} under $\cC$.
A collection $\cC$ is {\it triangularizable} 
if there is a chain of invariant subspaces for $\cC$ 
which is maximal as a subspace chain.

In many situations the existence of a non-trivial invariant subspace already 
implies triangularizability, as the Triangularization Lemma shows (see \cite{RR97} or \cite{RR00}). 
In order to recall it, some definitions are needed. 
Let $\cC$ be a collection of operators in $\BX$. 
If $\cM$ and $\cN$ are invariant subspaces under $\cC$ with $\cN \subset \cM$, then $\cC$
induces a collection $\hcC$ of quotients as follows: for
$T \in \cC$, the operator $\hT \in \hcC$ is defined on 
$\cM / \cN$ by
$$ \hT (x + \cN) \ = \ Tx + \cN \ . $$
Any such $\hcC$ is called a {\it set of quotients} of the collection $\cC$.
A property of collections of operators is said to be {\it inherited by
quotients} if every set of quotients of a collection having the property also has the same property.

\begin{lemma} (The Triangularization Lemma) Let $\cP$ be a property of collections of operators
that is inherited by quotients.  If every collection of operators (on a space of dimension greater than one)
which satisfies $\cP$ has a non-trivial invariant subspace,
then every collection satisfying $\cP$ is triangularizable.
\end{lemma}

A collection of operators is called an {\it (operator) semigroup} if it is closed under multiplication.
It is well-known fact that every commutative semigroup of matrices is triangularizable. 
In 1978 Laffey extended it as follows (see \cite{La78}).

\begin{theorem}
If $\cS$ is a semigroup of matrices such that the rank of $[S, T]$ is 
at most $1$ for all $\{S,T\} \subset \cS$, then $\cS$ is triangularizable.
\label{Laffey}
\end{theorem}

As a generalization of the preceding theorem, Radjavi and Rosenthal proved the following theorem (see \cite[Corollary 2]{RR97}).

\begin{theorem}
If $\cS$ is a semigroup contained in the Schatten class $C_p$ such that the rank of $[S, T]$ is at most $1$ for all $\{S,T\} \subset \cS$,
then $\cS$ is triangularizable.
\label{compact}
\end{theorem}

Applying the remarkable result of Turovskii \cite{Tu99} the last theorem can be easily extended to compact operators on $X$ (see \cite[Theorem 9.2.10]{RR00}).

Another infinite-dimensional extension of Laffey's result was shown by the group of authors in \cite{CD98}.
Recall that an operator $T$ on $X$ is said to be {\it algebraic} if there exists 
a nonzero complex polynomial $p$ such that $p(T) = 0$.

\begin{theorem}
Let $\cS$ be a semigroup of algebraic operators on $X$
such that the rank of $[S, T]$ is at most $1$ for all $\{S,T\} \subset \cS$. 
Then $\cS$ is triangularizable.
\label{algebraic}
\end{theorem}

In \cite{Dr00} we showed the following theorem related to the subject.

\begin{theorem}
Let $\cS$ be a non-commutative semigroup in $\BX$ generated by two elements.
If the rank of $[S, T]$ is at most $1$ for all $\{S,T\} \subset \cS$, then 
$\cS$ has a non-trivial hyperinvariant subspace.
\end{theorem}

The preceding theorem is an easy consequence of the following theorem \cite[Theorem 2.2]{Dr00}.

\begin{theorem}
Let $A$ and $B$ be operators on $X$ such that the rank of
$[A,B]$ is $1$ and the ranks of each of $[A^2,B]$, $[A,B^2]$, 
$[A^2,B^2]$, $[A B, B A]$ are at most $1$.  Then 
$\{A, B\}$ has a non-trivial hyperinvariant subspace.
\label{hyper}
\end{theorem}

It is clear from the proof of Theorem \ref{hyper} that the following is also true.

\begin{theorem}
Under the assumptions of Theorem \ref{hyper}, let $[A, B] = y \otimes \phi$ for a non-zero vector $y \in X$ 
and a non-zero functional $\phi \in X^*$. Then $\phi(y) = 0$, and so $[A, B]^2 = 0$.
\label{zero}
\end{theorem}

The proof of Theorem \ref{hyper} is essentially based upon the following
(easily proved) observation \cite[Lemma 2.1]{Dr00}.

\begin{lemma}
Let $y \in X$ and $\phi \in X^*$ be non-zero vectors.
Assume that $z \in X$ and $\psi \in X^*$ are such that 
the rank of the sum $y \otimes \phi + z \otimes \psi$ is at most $1$.
Then either $z$ is a multiple of $y$ or 
$\psi$ is a multiple of $\phi$.
\label{multiple}
\end{lemma}

\section{Results}

The key result of this paper is the following theorem. 

\begin{theorem}
Let $\cS$ be a semigroup in $\BX$ such that the rank of $[S, T]$ is at most $1$ for all $\{S,T\} \subset \cS$.
Suppose that $[A, B] = y \otimes \phi$ for some $\{A,B\} \subset \cS$ and for some non-zero vector $y \in X$ and 
non-zero functional $\phi \in X^*$. Then $\phi (C y) = 0$ for all $C \in \cS \cup \{I\}$.  
\label{main}
\end{theorem}

\begin{proof}
With no loss of generality we can assume that $\Cc \cS = \cS$, that is, $\cS$ 
is closed for scalar multiples of its members.
Since $\phi(y) = 0$ by Theorem \ref{zero}, we must prove that $\phi (C y) = 0$ for all $C \in \cS$. 
Assume on the contrary that $\phi (C y) \neq 0$ for some $C \in \cS$. 
As $\Cc \cS = \cS$, we can assume that $\phi (C y) = 1$.
Since $\phi(y) = 0$ and $(C^* \phi)(y) = \phi(C y) = 1$, neither $y$ is an eigenvector
of $C$ nor $\phi$ is an eigenvector of $C^*$. 
Therefore, it follows from Lemma \ref{multiple} that the rank of
$$ [AB,C] - [BA,C] = [[A,B],C] =  (y \otimes \phi)C - C (y \otimes \phi) =
y \otimes (C^*\phi) - (C y) \otimes \phi $$
is equal to $2$. Since the ranks of $[AB,C]$ and $[BA,C]$ are at most $1$, they are both equal to $1$.
In fact, we must have 
$$ [AB,C] = (\alpha y + \beta Cy) \otimes (\gamma \phi + \delta C^* \phi) $$
for some scalars $\alpha$, $\beta$, $\gamma$ and $\delta$.
We claim that $\delta \neq 0$. Assume on the contrary that $\delta = 0$.
Clearly, we can assume that $\gamma = 1$.
An application of Theorem \ref{zero} for the pair $\{AB, C\}$ yields 
$\phi(\alpha y + \beta Cy) = 0$, and so $\beta = 0$. 
Then the rank of 
$$ [BA, C] = [AB,C] - [[A,B],C] = y \otimes (\alpha \phi - C^*\phi) + (C y) \otimes \phi $$
is equal to $2$. This contradiction proves the claim.
Now, with no loss of generality we can assume that $\delta = 1$, and so 
$$ [AB,C] = (\alpha y + \beta Cy) \otimes (\gamma \phi + C^* \phi) \ \ \ \textrm{and} $$
$$ [BA, C] = [AB,C] - [[A,B],C] = y \otimes (\alpha \gamma \phi + (\alpha -1) C^*\phi) 
+ C y \otimes ((1 + \beta \gamma) \phi + \beta C^*\phi) . $$
Since $\rank [BA, C] \le 1$ and $C y$ is not a multiple of $y$, by Lemma \ref{multiple} there exists a scalar $\lambda$ such that 
$$ \alpha \gamma \phi + (\alpha -1) C^*\phi = \lambda ((1 + \beta \gamma) \phi + \beta C^*\phi) . $$
Since $C^*\phi$ and  $\phi$ are linearly independent, we have  
$$ \lambda (1 + \beta \gamma)) = \alpha \gamma \ \ \textrm{and} \ \ \lambda \beta = \alpha - 1 . $$
Eliminating $\alpha$ we obtain that $\lambda = \gamma$, and so $\alpha = 1 + \beta \gamma$.
An application of Theorem \ref{zero} for the pair $\{AB, C\}$ yields 
$(\gamma \phi + C^* \phi)(\alpha y + \beta Cy) = 0$, which implies that 
$$ 1 + 2 \beta \gamma + \beta \phi(C^2 y) = 0 . $$
It follows that $\beta \neq 0$, and we may define $k = 1/\beta + \gamma$.
Therefore, we have 
$$ [AB, C] = \beta (k y + Cy) \otimes (\gamma \phi + C^* \phi) , $$
$$ [BA, C] = \beta (\gamma y + Cy) \otimes (k \phi +  C^* \phi) $$
and 
$$ \phi(C^2 y) = - k - \gamma . $$
Now, the rank of 
$$  [AB, C^2] = [AB, C]C + C [AB, C] = $$
$$ = \beta (k y + Cy) \otimes (\gamma C^* \phi + (C^*)^2 \phi)  + 
\beta (k C y + C^2y) \otimes (\gamma \phi + C^* \phi)  $$
is at most $1$, and so, by Lemma \ref{multiple} again, one of the following two cases must occur. \\

\noindent{\it Case (I):} 
$k C y + C^2y = \lambda (k y + Cy)$ for some scalar $\lambda$.
Applying the functional $\phi$ to this equation we obtain that 
$\lambda = k + \phi(C^2y)$.  Since $\phi(C^2 y) = - k - \gamma$, we have  
$\lambda = - \gamma$, so that 
$$ C^2 y + (k + \gamma) C y + k \gamma y = 0  \ \ \ \textrm{or} \ \ \ (C + k) (C + \gamma) y  = 0 . $$
Now, decompose $X$ as $X = \lin \{y\} \oplus  \lin \{C y\} \oplus \ker (\phi) \cap \ker (C^* \phi)$.
Note that $\ker (\phi) = \lin \{y\} \oplus \ker (\phi) \cap \ker (C^* \phi)$.
With respect to this decomposition the operator $C$ has the matrix
$$ C = \left( \begin{array}{ccc}
  0  &  -k \gamma    & c_{1 3}  \\
  1  &  -k - \gamma  & c_{2 3}  \\
  0  &     0         & c_{3 3} 
  \end{array} \right) . $$
If $x \in  \ker (\phi) \cap \ker (C^* \phi)$, then $\phi (C x ) = 0$, so that 
$C x \in   \ker (\phi) = \lin \{y\} \oplus \ker (\phi) \cap \ker (C^* \phi)$. This implies that 
$c_{2 3} = 0$.
Since $[A, B] x = 0$ for all $x \in \ker (\phi)$ and $[A, B] Cy = y$, we have
$$ [A, B] = \left( \begin{array}{ccc}
  0  &  1  &  0  \\
  0  &  0  &  0  \\
  0  &  0  &  0  
  \end{array} \right) . $$
Since 
$$ [AB, C] y = \beta (k y + C y) $$
and 
$$ [AB, C] C y = - \beta k (k y + C y) , $$
we have 
\be [AB, C] = \beta \left( \begin{array}{ccc}
  k  &  -k^2 &  0  \\
  1  &   -k  &  0  \\
  0  &    0  &  0  
  \end{array} \right) = \beta \left( \begin{array}{ccc}
  k   \\
  1   \\
  0    
  \end{array} \right) 
  \left( \begin{array}{ccc}
  1  &   -k  &  0  
  \end{array} \right) . 
\label{AB_C} 
\ee  
Similarly, we obtain that 
$$ [BA, C] = \beta \left( \begin{array}{ccc}
  \gamma  &  -\gamma^2 &  0  \\
  1  &   -\gamma  &  0  \\
  0  &    0  &  0  
  \end{array} \right) = \beta \left( \begin{array}{ccc}
  \gamma   \\
  1   \\
  0    
  \end{array} \right) 
  \left( \begin{array}{ccc}
  1  &   -\gamma  &  0  
  \end{array} \right) . $$
Denoting $D = AB = (d_{i j})_{i, j = 1}^3$, we have 
$$ [AB, BA] = [A, B]D - D[A, B] = \left( \begin{array}{ccc}
  d_{2 1} & d_{2 2} - d_{1 1} & d_{2 3}  \\
    0     &  - d_{2 1}  &  0  \\
    0     &  - d_{3 1}  &  0  
  \end{array} \right) . $$
Since $\rank [AB, BA] \le 1$, we conclude that $d_{2 1} = 0$, and either     
$d_{2 3} = 0$ or $d_{3 1} = 0$. Therefore, we must consider two subcases. \\

\noindent{\it Subcase (Ia):} $d_{2 3} = 0$. Then 
$$ [A B, C] = [A B - d_{1 1} I, C] = $$
$$ = \left( \begin{array}{ccc}
 d_{1 2} - c_{1 3} d_{3 1} & -(k + \gamma)d_{1 2} + k \gamma (d_{2 2} - d_{1 1}) - c_{1 3} d_{3 2} & d_{1 3} c_{3 3} - c_{1 3} (d_{3 3} - d_{1 1}) \\
 d_{2 2} - d_{1 1} &  - d_{1 2} & - d_{1 3} \\
 d_{3 2} - c_{3 3} d_{3 1} & - k \gamma d_{3 1} - (k + \gamma)d_{3 2} - c_{3 3} d_{3 2} & 
                                               d_{3 1} c_{1 3} + d_{3 3} c_{3 3} -  c_{3 3} d_{3 3}
\end{array} \right) . $$
Comparing with (\ref{AB_C}) and simplifying we get
$$ d_{1 3} = 0 , \ d_{1 2} = \beta k , \ d_{2 2} = d_{1 1} + \beta , $$
\be c_{13} d_{31} = 0 , \ c_{13} d_{32}= 0 , 
\ d_{32} = c_{33} d_{31} . 
\label{cd} \ee
We thus have
$$ D = A B =  \left( \begin{array}{ccc}
  d_{1 1} & \beta k &  0  \\
    0     & d_{1 1} + \beta &  0  \\
  d_{3 1} &  d_{3 2}  &  d_{3 3}   
  \end{array} \right)  \ \ \ \textrm{and} $$
$$ B A =  A B - [A,B] = \left( \begin{array}{ccc}
  d_{1 1} & \beta \gamma &  0  \\
    0     & d_{1 1} + \beta &  0  \\
  d_{3 1} &  d_{3 2}  &  d_{3 3}   
  \end{array} \right) . $$
Now compute the commutator
$$ [AB BA, C] = [A B, C]B A + A B[B A, C] = $$
$$ = \beta \left( \begin{array}{ccc}
  k   \\
  1   \\
  0    
  \end{array} \right) 
  \left( \begin{array}{ccc}
  1  &   -k  &  0  
  \end{array} \right)
  \left( \begin{array}{ccc}
  d_{1 1} & \beta \gamma &  0  \\
    0     & d_{1 1} + \beta &  0  \\
  d_{3 1} &  d_{3 2}  &  d_{3 3}   
  \end{array} \right) + $$
$$ +  \left( \begin{array}{ccc}
  d_{1 1} & \beta k &  0  \\
    0     & d_{1 1} + \beta &  0  \\
  d_{3 1} &  d_{3 2}  &  d_{3 3}   
  \end{array} \right)    
\beta \left( \begin{array}{ccc}
  \gamma  \\
  1   \\
  0    
  \end{array} \right) 
  \left( \begin{array}{ccc}
  1  &   -\gamma  &  0  
  \end{array} \right) = $$
$$ = \beta \left( \begin{array}{ccc}
  k   \\
  1   \\
  0    
  \end{array} \right) 
  \left( \begin{array}{ccc}
  d_{1 1}  &   - k d_{1 1} - 1 &  0  
  \end{array} \right)
+ 
  \beta \left( \begin{array}{ccc}
  d_{1 1} \gamma + \beta k  \\
  d_{1 1} + \beta  \\
  d_{3 1} \gamma + d_{3 2}   
  \end{array} \right)    
  \left( \begin{array}{ccc}
  1  &   - \gamma  &  0  
  \end{array} \right) . $$
Since its rank is at most $1$, two subsubcases are possible by Lemma \ref{multiple}. \\

\noindent{\it Subsubcase (Ia1):}
$(d_{1 1} \gamma + \beta k, d_{1 1} + \beta, d_{3 1} \gamma + d_{3 2}) = \mu (k, 1, 0)$   
for some scalar $\mu$. Eliminating $\mu$ we obtain that $d_{11} = 0$ and 
\be d_{3 1} \gamma + d_{3 2} = 0 .  \label{d31d32} \ee
We now compute the commutator
$$ [BA AB, C] = [BA, C]A B + B A [A B, C] = $$
$$ = \beta \left( \begin{array}{ccc}
  \gamma   \\
  1   \\
  0    
  \end{array} \right) 
  \left( \begin{array}{ccc}
  1  &   -\gamma  &  0  
  \end{array} \right)
  \left( \begin{array}{ccc}
    0   & \beta k &  0  \\
    0   & \beta &  0  \\
  d_{3 1} &  d_{3 2}  &  d_{3 3}   
  \end{array} \right) + $$
$$ +  \left( \begin{array}{ccc}
    0   & \beta \gamma &  0  \\
    0   & \beta &  0  \\
  d_{3 1} &  d_{3 2}  &  d_{3 3}   
  \end{array} \right)    
\beta \left( \begin{array}{ccc}
  k  \\
  1   \\
  0    
  \end{array} \right) 
  \left( \begin{array}{ccc}
  1  &  - k  &  0  
  \end{array} \right) = $$
$$ = \beta \left( \begin{array}{ccc}
  \gamma   \\
  1   \\
  0    
  \end{array} \right) 
  \left( \begin{array}{ccc}
  0  &   1  &   0  
  \end{array} \right)
+ 
  \beta \left( \begin{array}{ccc}
  \beta \gamma \\
  \beta  \\
  d_{3 1} k + d_{3 2}   
  \end{array} \right)    
  \left( \begin{array}{ccc}
  1  &   - k &  0  
  \end{array} \right) . $$
Since its rank is at most $1$, we conclude that $d_{3 1} k + d_{3 2} = 0$. 
Using (\ref{d31d32}) we obtain that $d_{3 1} = 0$ and $d_{3 2} = 0$ as $k \neq \gamma$. 
Therefore, the two-dimensional subspace $\lin \{y, C y\}$ is invariant under $\{C, AB, BA \}$. 
Let $\tilde{C}$, $\tilde{D}$ and $\tilde{E}$ denote the restrictions to this subspace of 
$C$, $D = AB$ and $E=BA$, respectively. Then 
$$ \tilde{C} = \left( \begin{array}{ccc}
  0  &  -k \gamma  \\
  1  &  -k - \gamma  
  \end{array} \right) , $$
$$ \tilde{D} = \beta \left( \begin{array}{ccc}
  k   \\
  1       
  \end{array} \right) 
  \left( \begin{array}{ccc}
  0 &  1    
  \end{array} \right)  \ \ \ \textrm{and} \ \ \ 
\tilde{E} = \beta \left( \begin{array}{ccc}
  \gamma   \\
  1       
  \end{array} \right) 
  \left( \begin{array}{ccc}
  0 &  1    
  \end{array} \right) . $$
Now, we could apply Theorem \ref{Laffey} to obtain 
triangularizability of the semigroup generated by $\tilde{C}$, $\tilde{D}$ and $\tilde{E}$ 
which would contradict the fact that the commutator 
$$ [\tilde{D}-\tilde{E}, \tilde{C}] = 
  \left( \begin{array}{ccc}
  1  &  -k - \gamma  \\
  0  &  -1  
  \end{array} \right) $$
is not nilpotent. However, there exists a quick way to get a contradiction directly. Since  
$$ \tilde{D} \tilde{C} = 
\beta \left( \begin{array}{ccc}
  k   \\
  1       
  \end{array} \right) 
  \left( \begin{array}{ccc}
  1 & -k-\gamma      
  \end{array} \right) \ \ \ \textrm{and} \ \ \ 
\tilde{E} \tilde{C} = 
\beta \left( \begin{array}{ccc}
  \gamma   \\
  1       
  \end{array} \right) 
  \left( \begin{array}{ccc}
  1 & -k-\gamma      
  \end{array} \right) , $$
we have  
$$ [ \tilde{D} \tilde{C},  \tilde{E}] = 
   \beta^2 \left\{ \left( \begin{array}{ccc}
  k  \\
  1       
  \end{array} \right) 
  \left( \begin{array}{ccc}
  0 &  -k    
  \end{array} \right) 
+ 
  \left( \begin{array}{ccc}
  \gamma  \\
  1       
  \end{array} \right) 
  \left( \begin{array}{ccc}
  1 &  k+\gamma     
  \end{array} \right) \right\} \ \ \ \textrm{and} $$
$$ [ \tilde{E} \tilde{C},  \tilde{D}] = 
   \beta^2 \left\{ \left( \begin{array}{ccc}
  \gamma  \\
  1       
  \end{array} \right) 
  \left( \begin{array}{ccc}
  0 &  -\gamma  
  \end{array} \right) 
+
  \left( \begin{array}{ccc}
  k  \\
  1       
  \end{array} \right) 
  \left( \begin{array}{ccc}
  1 &  k+\gamma     
  \end{array} \right) \right\} . $$
Since the ranks of both commutators are at most $1$, we conclude first that $k = 0$ and then $\gamma = 0$.
This is not possible, and so the subsubcase (Ia1) is finished. \\

\noindent{\it Subsubcase (Ia2):}
$(d_{1 1}, - k d_{1 1} - 1, 0) = \mu (1, - \gamma, 0)$ 
for some scalar $\mu$. It follows that $d_{11} = \mu = -\beta$.
In this case we compute   
$$ [C AB, BA] = C [A B, B A] - [B A, C] AB  = $$ 
$$ = \left( \begin{array}{ccc}
  0  &  -k \gamma    & c_{1 3}  \\
  1  &  -k - \gamma  & 0        \\
  0  &     0         & c_{3 3} 
  \end{array} \right)  
  \left( \begin{array}{ccc}
  \beta  \\
  0      \\
  - d_{31}      
  \end{array} \right) 
  \left( \begin{array}{ccc}
  0 & 1 & 0      
  \end{array} \right) - $$
$$ - \beta \left( \begin{array}{ccc}
  \gamma   \\
  1   \\
  0    
  \end{array} \right) 
  \left( \begin{array}{ccc}
  1  &   -\gamma  &  0  
  \end{array} \right) 
  \left( \begin{array}{ccc}
  - \beta & \beta k &  0  \\
    0     & 0   &    0  \\
  d_{3 1} &  d_{3 2}  &  d_{3 3}   
  \end{array} \right) = $$
$$  =
\left( \begin{array}{ccc}
   0  \\
   \beta  \\
  - d_{32}      
  \end{array} \right) 
  \left( \begin{array}{ccc}
  0 & 1 & 0      
  \end{array} \right) - 
\beta \left( \begin{array}{ccc}
  \gamma   \\
  1   \\
  0    
  \end{array} \right) 
  \left( \begin{array}{ccc}
  -\beta  &  \beta k  &  0  
  \end{array} \right), $$
 where we have used two equalities from (\ref{cd}).
Since the rank of this commutator is at most $1$,  
we obtain that $\gamma = 0$ and $d_{32}=0$.
Similarly, we obtain that
$$ [C B A, A B] = C [B A, A B] - [AB, C] BA  = $$ 
$$ = \left( \begin{array}{ccc}
  0  &  0    & c_{1 3}  \\
  1  &  -k   & 0        \\
  0  &  0    & c_{3 3} 
  \end{array} \right)  
  \left( \begin{array}{ccc}
  - \beta  \\
  0      \\
  d_{31}      
  \end{array} \right) 
  \left( \begin{array}{ccc}
  0 & 1 & 0      
  \end{array} \right) - $$
$$ - \beta \left( \begin{array}{ccc}
  k  \\
  1   \\
  0    
  \end{array} \right) 
  \left( \begin{array}{ccc}
  1  &   - k  &  0  
  \end{array} \right) 
  \left( \begin{array}{ccc}
  - \beta & 0 &  0  \\
    0     & 0  &    0  \\
  d_{3 1} & 0  &  d_{3 3}   
  \end{array} \right) = $$
$$  = 
\left( \begin{array}{ccc}
   0  \\
 -\beta  \\
   0       
  \end{array} \right) 
  \left( \begin{array}{ccc}
  0 & 1 & 0      
  \end{array} \right) - 
\beta \left( \begin{array}{ccc}
  k   \\
  1   \\
  0    
  \end{array} \right) 
  \left( \begin{array}{ccc}
  \beta  &  0  &  0  
  \end{array} \right) .$$
It follows that $k = 0$, so that $\gamma = k = 0$.
This contradiction completes the proof in this subsubcase,
and so Subcase (Ia) as well. \\

\noindent{\it Subcase (Ib):} $d_{31} = 0$. Then 
$$ [A B, C] = \left( \begin{array}{ccc}
 d_{11} &  d_{12} & d_{13} \\
   0    &  d_{22} & d_{23} \\
   0    &  d_{32} & d_{33} 
 \end{array} \right) 
\left( \begin{array}{ccc}
  0  &  -k \gamma    & c_{1 3}  \\
  1  &  -k - \gamma  & 0        \\
  0  &     0         & c_{3 3} 
  \end{array} \right) - $$
$$ - 
\left( \begin{array}{ccc}
  0  &  -k \gamma    & c_{1 3}  \\
  1  &  -k - \gamma  & 0        \\
  0  &     0         & c_{3 3} 
  \end{array} \right)   
  \left( \begin{array}{ccc}
 d_{11} &  d_{12} & d_{13} \\
   0    &  d_{22} & d_{23} \\
   0    &  d_{32} & d_{33} 
 \end{array} \right) = 
 \left( \begin{array}{ccc}
 d_{12} &  *  &  * \\
 d_{22} - d_{11}  &  *  &  * \\
 d_{32} &  *  &  * 
 \end{array} \right) . $$
Comparing with (\ref{AB_C}) we obtain $d_{12} = \beta k$, $d_{22} = d_{11} + \beta$ and $d_{32} = 0$. 
Therefore, the two-dimensional subspace $\lin \{y, C y\}$ is invariant under $\{C, AB, BA \}$,
and so we can proceed as in Subsubcase (Ia1) to get a contradiction.
Denote by $\tilde{C}$, $\tilde{D}$ and $\tilde{E}$ the restrictions to this invariant subspace of 
$C$, $D = AB$ and $E=BA$, respectively. Then 
$$ \tilde{C} = \left( \begin{array}{ccc}
  0  &  -k \gamma  \\
  1  &  -k - \gamma  
  \end{array} \right) , $$
$$ \tilde{D} = \left( \begin{array}{ccc}
  d_{11} & \beta k   \\
  0      & d_{11} + \beta 
  \end{array} \right) 
  \ \ \ \textrm{and} \ \ \ 
\tilde{E} = \left( \begin{array}{ccc}
  d_{11} & \beta \gamma   \\
  0      & d_{11} + \beta 
  \end{array} \right) . $$
By Theorem \ref{Laffey}, the semigroup generated by $\{ \tilde{C}, \tilde{D}, \tilde{E}\}$
is triangularizable. However, this implies a contradiction, as the commutator 
$$ [\tilde{D}-\tilde{E}, \tilde{C}] = 
  \left( \begin{array}{ccc}
  1  &  -k - \gamma  \\
  0  &  -1  
  \end{array} \right) $$
is not nilpotent. 
This concludes this subcase, and so Case (I) as well. \\

\noindent{\it Case (II):} 
$\gamma C^* \phi + (C^*)^2 \phi = \lambda (\gamma \phi + C^* \phi)$ for some scalar $\lambda$. 
Computing both sides of this equality at the vector $y$ yields
$\lambda = \gamma + \phi (C^2 y)$. Since $\phi (C^2 y) = - k - \gamma$, we obtain that 
$\lambda = -k$. Therefore, we have 
$$ (C^*)^2 \phi + (k+\gamma) C^* \phi + k \gamma \phi = 0 . $$
Thus, the semigroup $\cS^* = \{S^* : S \in \cS\}$ satisfies the same conditions as 
the semigroup $\cS$ in Case (I). Indeed, if $x \mapsto F_x$ denotes the isometric embedding of $X$ to $X^{**}$,
then we have $[B^*, A^*] = \phi \otimes F_y$, $F_y(\phi) = \phi(y) = 0$ and 
$F_y(C^* \phi) = (C^* \phi)(y) = \phi(C y) \neq 0$. 
So, we can use Case (I) to obtain a contradiction. This completes the proof of the theorem.
\end{proof}

\begin{corollary}
If $\cS$ is a non-commutative semigroup in $\BX$ such that the rank of $[S, T]$ is at most $1$ for all $\{S,T\} \subset \cS$,
then $\cS$ has a non-trivial invariant subspace.
\label{non-comm}
\end{corollary}

\begin{proof}
Since $\cS$ is not commutative, we have $[A, B] = y \otimes \phi$ for some $\{A,B\} \subset \cS$ and for some non-zero vector 
$y \in X$ and non-zero functional $\phi \in X^*$. 
Then the closed linear span of the set $\{S y : S \in \cS \cup \{I\}\}$ 
is a non-zero subspace invariant under $\cS$. Since $\phi (S y) = 0$ for all $S \in \cS \cup \{I\}$ 
by Theorem \ref{main}, it is contained in $\ker (\phi)$, and so it is non-trivial.
\end{proof}

The preceding corollary does not hold for commutative semigroups, since $\cS$ may be generated by an operator
without non-trivial invariant subspaces. In order to cover the commutative case as well, 
the following property of operators was introduced in \cite{Dr00}.
Let $\cR$ be the property of operators on Banach spaces of dimension
at least $2$ such that:

(a) $\cR$ is inherited by quotients,

(b) each commutative semigroup of operators with the property $\cR$
has a non-trivial invariant subspace.

Since every non-zero compact operator on an infinite-dimensional Banach space has a non-trivial hyperinvariant 
subspace by the famous Lomonosov's result \cite{Lo73}, the property of being a compact operator is an example
of such property $\cR$. Another example is the property of being an algebraic operator, as 
eigenspaces of an algebraic operator are hyperinvariant subspaces. \\

The main result of the paper is the following theorem.

\begin{theorem}
Let $\cS$ be a semigroup of operators on $X$ with the property $\cR$.  
If the rank of $[S, T]$ is at most $1$ for all $\{S, T\} \subset \cS$, 
then $\cS$ is triangularizable.
\label{triang}
\end{theorem}

\begin{proof}
By the Triangularization Lemma, it suffices to show that $\cS$ has 
a nontrivial invariant subspace. This is true by Corollary \ref{non-comm}  
if $\cS$ is not commutative. Otherwise, this holds by the condition (b)
of property $\cR$.
\end{proof}

Compact and algebraic operators are special cases of polynomially compact operators. Recall that an operator $T$ on $X$ 
is said to be {\it polynomially compact} if there exists a nonzero complex polynomial $p$ such that the operator $p(T)$ 
is compact. Triangularizability of collections of polynomially compact operators was studied by 
Konvalinka in \cite{Ko05}. Since the Lomonosov's result is strong enough to give 
non-trivial hyperinvariant subspaces of polynomially compact operators (that are not algebraic), 
the property of being a polynomially compact operator is also an example of property $\cR$. So, we have the following 
corollary of Theorem \ref{triang} that extends both Theorems \ref{compact} and \ref{algebraic}.

\begin{corollary}
Let $\cS$ be a semigroup of polynomially compact operators on $X$.  
If the rank of $[S, T]$ is at most $1$ for all $\{S, T\} \subset \cS$, 
then $\cS$ is triangularizable.
\end{corollary}


\begin{thebibliography}{9999}

\bibitem{CD98} G. Cigler, R. Drnov\v sek, D. Kokol-Bukov\v sek, T. Laffey, 
M. Omladi\v c, H. Radjavi, P. Rosenthal, 
{\it Invariant subspaces for semigroups of algebraic operators}, 
J. Funct. Anal. {\bf 160} (1998), 452--465.

\bibitem{Dr00} R. Drnov\v sek, {\it Hyperinvariant subspaces for operator semigroups with commutators of rank at most one}, 
Houston J. Math.  {\bf 26} (2000),  no. 3, 543--548.

\bibitem{Ko05} M. Konvalinka, {\it Triangularizability of polynomially compact operators}, 
Integr. Equ. Oper. Theory {\bf 52} (2005), no. 2, 271--284. 

\bibitem{La78} T. J. Laffey, {\it Simultaneous triangularization of matrices
- low rank cases and the nonderogatory case}, Lin. Mult. Alg. {\bf 6} (1978), 269--305.

\bibitem{Lo73} V. I. Lomonosov, {\it Invariant subspaces of the family of operators that commute with a completely continuous operator},
Funkcional. Anal. i Prilo\v zen.  {\bf 7} (1973), no. 3, 55--56. 

\bibitem{RR97} H. Radjavi and P. Rosenthal, {\it From local to global
triangularization}, J. Funct. Anal. {\bf 147} (1997), 443--456.

\bibitem{RR00} H. Radjavi and P. Rosenthal, {\it Simultaneous triangularization}, 
Universitext, Springer-Verlag, New York, 2000.

\bibitem{Tu99} Yu. V. Turovskii,
{\it Volterra semigroups have invariant subspaces},
J. Funct. Anal. {\bf 162} (1999), 313--322.


\end{thebibliography}
\end{document}